\documentclass[onecolumn, 11 pt]{ieeeconf}  % Use a4paper instead of letterpaper if you need to change size of the paper

\IEEEoverridecommandlockouts                              % This command is only needed if 
\usepackage{amsmath,amssymb}
\usepackage{mathptmx}
\usepackage{graphics,graphicx,epsfig}
\usepackage{cite}
\usepackage[usenames,dvipsnames,svgnames,table]{xcolor}
\usepackage{dsfont}

\newcommand{\R}{\mathbb{R}} 
\newcommand{\E}{\mathbb{E}}

\usepackage{flushend}

\usepackage{comment}

\newif\ifshow
\showtrue
\ifshow
\newenvironment{scratch}
{
  \color{ForestGreen}
}
{
}
\else
\excludecomment{scratch}
\fi

\usepackage{cite}

\newtheorem{remark}{\bfseries Remark}

\title{\LARGE \bf
Driving an Ornstein--Uhlenbeck Process to Desired \\ First-Passage Time Statistics}
\author{Khem Raj Ghusinga$^{1}$, Vaibhav Srivastava$^{2}$, and Abhyudai Singh$^{1}$ % <-this % stops a space
\thanks{$^{1}$ Khem Raj Ghusinga and Abhyudai Singh are with the Department of Electrical and Computer Engineering, University of Delaware, Newark, DE, USA.
        {\tt\small \{khem,absingh\}@udel.edu}}%
\thanks{$^{2}$ Vaibhav Srivastava is with the Department of Electrical and Computer Engineering, Michigan State University, East Lansing, MI, USA.
        {\tt\small vaibhav@engr.msu.edu}}%
}

\begin{document}

\maketitle
\thispagestyle{empty}
\pagestyle{empty}

%%%%%%%%%%%%%%%%%%%%%%%%%%%%%%%%%%%%%%%%%%%%%%%%%%%%%%%%%%%%%%%%%%%%%%%%%%%%%%%%
\begin{abstract}
First-passage time (FPT) of an Ornstein-Uhlenbeck (OU) process is of immense interest in a variety of contexts. This paper considers an OU process with two boundaries, one of which is absorbing while the other one could be either reflecting or absorbing, and studies the control strategies that can lead to desired FPT moments. Our analysis shows that the FPT distribution of an OU process is scale invariant with respect to the drift parameter, i.e., the drift parameter just controls the mean FPT and doesn't affect the shape of the distribution. This allows to independently control the mean and coefficient of variation (CV) of the FPT. We show that that increasing the threshold may increase or decrease CV of the FPT, depending upon whether or not one of the threshold is reflecting. 
We also explore the effect of control parameters on the FPT distribution, and find parameters that minimize the distance between the FPT distribution and a desired distribution. 
%We also find optimal drift and diffusion parameters that yield a FPT distribution close to a desired distribution. 
\end{abstract}

%%%%%%%%%%%%%%%%%%%%%%%%%%%%%%%%%%%%%%%%%%%%%%%%%%%%%%%%%%%%%%%%%%%%%%%%%%%%%%%%
\section{INTRODUCTION}

The first passage time (FPT) is the earliest time at which a trajectory of a stochastic process initially inside a bounded region leaves the region. The FPTs are extensively used across disciplines, including neuroscience~\cite{bogacz2006physics}, biology~\cite{HCT-FYMW:00, GV-BF-JC-AD-ERD:08, SS-DJS-CRM:14}, finance~\cite{redner2001guide}, ecology~\cite{OB-MC-MM-PHS-RV:05}, engineering~\cite{WS-MZW-LAS:12}, statistical physics~\cite{Farkas2001}, finance~\cite{Lin1998}, and health science~\cite{Horrocks2004} to model several interesting phenomena. For example, the FPT of diffusion processes is used to model human decision-making~\cite{bogacz2006physics}, animal foraging~\cite{OB-MC-MM-PHS-RV:05}, financial markets~\cite{Lin1998}, and clock synchronization~\cite{WS-MZW-LAS:12}. 

In this paper, we study control of the FPT statistics of an Ornstein-Uhlenbeck (OU) process between two fixed boundaries; one of which is absorbing, and the other can be absorbing or reflecting. An OU process belongs to the class of diffusion processes and is a generalization of drift-diffusion process. The OU process is also a continuum approximation to several discrete time Markov models. For biological phenomena modeled by FPTs, the analysis in this paper can provide insights into the mechanisms these systems employ to cope with uncertainty and ensure resilient performance. For example, how attention and memory is modulated in human decision-making, or how a gene's expression is regulated to control timing of its response, or how animals regulate their foraging activity. For engineered systems, these analysis may provide insights into optimal control laws that delay an undesired event such as epidemic outbreak, or optimal control laws that achieve a desired distribution for time to certain event such as adoption of a product by certain fraction of population.

The problem of steering a linear stochastic system to a desired final distribution has been studied~\cite{chen2016optimal}. However, computing and controlling FPT distribution is significantly more complicated than controlling the evolution of trajectories without boundaries. Indeed, the Fokker-Plank equation for the OU process is nonlinear and has limited tractability~\cite{gardiner2004handbook}. Control of FPT distribution for OU process has been studied in~\cite{sacerdote2003threshold}, wherein the boundary of the region is controlled to steer the FPT distribution to a Gamma distribution. Loosely speaking, this problem can be thought of as a boundary control of a PDE~\cite{krstic2008boundary}, where underlying PDE is the Fokker-Planck equation. 

In our recent work~\cite{Ghusinga17}, similar problems were explored in the context of gene expression. Therein the stochastic process is a continuous-time discrete-state  process defined on positive integers, with a reflecting boundary at $0$ and a fixed absorbing boundary. The results showed that the best strategy to minimize the coefficient of variation (CV) of FPT for a fixed mean FPT is a constant rate of production (forward hopping) and no decay (backward hopping). These results interpreted in continuum limit would mean that the optimal stochastic process (within OU processes) for minimizing the CV of FPT for a given mean is the drift-diffusion process. In other words, the optimal control is a feedforward controller and requires no state feedback. In this paper, we explore this control problem in more detail.

%\com{Motivation for absorbing, and reflected boundaries}

%\com{Related Work} % just trying to create outline of the introduction. we would comment this line later.

%\com{Contributions}
Although the FPT properties of OU processes have been extensively studied in the literature~\cite{BorodinSalminen12}, a control theoretic analysis of how the process can be steered to some desired FPT statistics is lacking. This paper provides an unified approach based on characteristic functions to find control parameters that lead to desired FPT distributions. The approach is analytically and numerically tractable and provides important insights into the FPT behavior of the OU process. The major contributions of this paper are threefold.
First, we show that the FPT distribution for OU process is scale-invariant with respect to drift parameter, which facilitates independent tuning of the mean and CV of the FPT. Second, using the characteristic function of the FPT, we explore the space of control parameters to understand the variation of the FPT statistics with these parameters. Third, we determine optimal control parameters to achieve desired FPT statistics.

%\com{Organization}

The paper is organized as follows. Section~\ref{sec:problem}
introduces the problem. Section~\ref{sec:charfun} presents background results on the characteristic function for the FPT of an OU process. Section~\ref{sec:OptMom}
uses the characteristic function to find properties of moments of the FPT, and optimal parameters that lead to  desired moments. A more general control problem that explores the parameter space to reach a desired FPT distribution is studied in Section~\ref{sec:OptDist}. Finally, conclusions and
future work are discussed in Section~\ref{sec:conclusion}.

\section{PROBLEM DESCRIPTION}\label{sec:problem}
Consider an  OU process defined by the following stochastic differential equation
\begin{equation} \label{eq:ou}
dx=-\theta x dt+\sigma \sqrt{\theta} dw_t.
\end{equation}
Here $x$ is the state, $\theta \in \R_{\ge 0}$, and $\sigma \in \R_{\ge 0}$ are parameters, and $dw_t$ are i.i.d. Wiener increments. We will refer to $\theta$ as the drift and $\sigma$ as the relative noise strength. Let $a$ and $b$ denote two thresholds such that $a<b$. The FPT, $\tau$, for $x(t)$ to cross either of these thresholds is mathematically defined as
\begin{equation} \label{eq:fpt2thresholds}
\tau=\inf\{t: x(t) \notin (a,b) | x(0)=x_0 \in (a,b)\}.
\end{equation}

% \begin{figure}[h]
% \centering
% \includegraphics[width=\linewidth]{fig1.pdf}
% \caption{FPT for an OU process with two absorbing thresholds.}
% \label{fig:OUSchematic}
% \end{figure}

Our aim is to investigate optimal drift $\theta$ and relative noise strength $\sigma$ that lead to desired FPT moments. Such problems could be of relevance in many contexts wherein a desired mean FPT and at least a tolerable CV is required. This problem can be generalized further by demanding a FPT distribution that is as close to a desired distribution as it could be. The thresholds $a$ and $b$ could be both absorbing, or one absorbing and the other reflecting.  
%This problem can be casted in different variants, such as having one of the thresholds to be reflecting.

To handle these problems in an unified manner, we propose to use the characteristic function of the OU process. Not only the moments can be easily computed from the characteristic function, but it also provides a useful way to characterize the distance between two probability distribution functions. To be more specific, the characteristic function, $\psi_\tau(\alpha)$, of the FPT, $\tau$, is defined as
\begin{equation}
\psi_\tau(\alpha)=\E \left[e^{i \alpha \tau}\right], \quad \alpha \in \R.
\end{equation}
A $m$-th order moment $\E[\tau^m]$ can be computed as
\begin{equation}
\E[\tau^m]=i^{-m}\left[\frac{d^m}{d\alpha^m}\psi_\tau(\alpha)\right]_{\alpha=0}.
\end{equation}
Furthermore, the following result due to Parseval--Plancherel provides a metric to quantify the difference between two probability density functions in terms of their characteristic functions:
Let $f_{\tau}(t)$ denote the probability density function of the FPT, and $f_d(t)$ be a desired probability distribution function. The distance between these functions can be quantified in terms of their characteristic functions $\psi_\tau(\alpha)$ and $\psi_d(\alpha)$ as
\begin{equation}\label{eq:distmetric}
\int_{0}^{\infty}\left[f_\tau(t)-f_d(t)\right]^2 dt = \frac{1}{2\pi}\int_{-\infty}^{\infty} \left \lvert \psi_{\tau}(\alpha)-\psi_{d}(\alpha)\right \rvert^2 d\alpha,
\end{equation}
provided that the integrals exist \cite{UshakovCF99}.

%%%%%%%%%%%%%%%%%%%%%%%%%%%%%%%%%%%%%%%%%%%%%%%%%%%%%%%%%%%%%%%%%%%%%%%%%%%%%%%%
\section{BACKGROUND RESULTS ON FPT OF AN OU PROCESS}\label{sec:charfun}
In this section, we provide background results on FPT of the OU process \eqref{eq:ou}. For completeness, we provide detailed computation of the characteristic function using standard tools from the theory of stochastic processes (see, \cite{BorodinSalminen12,GihmanSkorohodSDE72, CoxMiller77,gardiner2004handbook,redner2001guide}).
We consider two thresholds at $a$ and $b$, both of which could be absorbing or one of them could be reflecting.

\subsection{When both thresholds are absorbing}
To derive the characteristic function, $\psi_{\tau}(\alpha)$, for the OU process in \eqref{eq:ou}, we define $g(y)$ as
\begin{equation}
g(y)=\E \left[e^{i \alpha \tau(y)}\right],
\end{equation}
where $y$ represents an initial condition, and $\tau(y)$ denotes the FPT starting from an initial condition $y$. Note that the characteristic function is related with $g(y)$ as $\psi_{\tau}=g(x_0)$. The computation of $g(y)$ using first principles is discussed below.

Consider the evolution of the OU process starting from $y$ in an infinitesimal time interval $h$. Denote $x_h=x(h)=y-\theta y h+\sigma \sqrt{\theta} w_h$. It follows that 
\begin{subequations}
\begin{align}
g(y) &= \E_{x_h} \E_{\tau(x_h)}\left[e^{i \alpha(h+\tau(x_h))}\right] \\
&= e^{i \alpha h}\E_{x_h}\left[g(x_h)\right] \\
&= e^{i \alpha h} \left(g(y) - \theta y h \frac{dg}{dy} + \frac{1}{2} \sigma^2 \theta h \frac{d^2g}{dy^2}  \right) + O(h^2).
\end{align}
\end{subequations}
Taking the limit $h \to 0$ results in
\begin{equation}\label{eq:gODE}
\frac{1}{2}\sigma^2 \theta \frac{d^2g(y)}{dy^2} -\theta y \frac{dg(y)}{dy} + i\alpha g(y)=0.
\end{equation}

We are interested in the solution to the above differential equation which can be obtained using the series method. Let 
\begin{equation}\label{eq:series}
g(y)=\sum_{n=0}^{\infty}c_n y^n.
\end{equation}
Plugging this in \eqref{eq:gODE} results in
\begin{equation}\label{eq:Gdiffeqn}
\frac{1}{2}\sigma^2 \theta \sum_{n=0}^{\infty}(n+2)(n+1)c_{n+2}y^n -\theta \sum_{n=0}^{\infty}n c_{n} y^n + i\alpha \sum_{n=0}^{\infty}c_n y^n=0.
\end{equation}
It is straightforward to see that \eqref{eq:Gdiffeqn} results in the following recursive relation in the coefficients
\begin{equation}\label{eq:coeffs}
c_{n+2}=\frac{2(-i\alpha + n \theta)}{\sigma^2 \theta (n+2)(n+1)} c_n.
\end{equation}
The above recursion yields the following solution
\begin{subequations}\label{eq:coeffOddEven}
\begin{align}
& c_n=\frac{2^{n} \Gamma \left[\frac{1}{2}\left(n-\frac{i \alpha }{\theta}\right)\right]c_0}{\Gamma \left[-\frac{i \alpha }{2 \theta }\right]\sigma^n n!}, \quad n=0, 2, 4, \ldots \\
& c_n=\frac{2^{n-1} \Gamma \left[\frac{1}{2}\left(n-\frac{i \alpha }{\theta}\right)\right]c_1}{\Gamma \left[-\frac{i \alpha + \theta }{2 \theta }\right]\sigma^{n-1} n!}, \quad n=1, 3, 5, \ldots,
\end{align}
\end{subequations}
where $\Gamma$ is the Gamma function. A general solution to \eqref{eq:gODE} can be given by \eqref{eq:series}, with the coefficients given by \eqref{eq:coeffOddEven}. Simplifying the series in \eqref{eq:series} via symbolic manipulation in Mathematica yields 
\begin{equation}\label{eq:gysol}
g(y)=c_0 \, _1F_1\left(-\frac{i \alpha }{2 \theta },\frac{1}{2},\frac{y^2}{\sigma ^2}\right)+c_1 y \, _1F_1\left(\frac{\theta -i \alpha }{2 \theta }, \frac{3}{2},\frac{y^2}{\sigma ^2}\right),
\end{equation}
where $_1F_1$ represents the Kummer's confluent hypergeometric function.

The solution in \eqref{eq:gysol} consists of two unknown coefficients $c_0$ and $c_1$ which can be computed using the boundary conditions. When both thresholds $a$ and $b$ are absorbing, the boundary conditions are given by $g(a)=1$ and $g(b)=1$. Using these boundary values, $c_0$ and $c_1$ can be determined by solving
\begin{subequations}
\begin{align}
&c_0 \, _1F_1\left(-\frac{i \alpha }{2 \theta },\frac{1}{2},\frac{a^2}{\sigma ^2}\right)+c_1 a \, _1F_1\left(\frac{\theta -i \alpha }{2 \theta }, \frac{3}{2},\frac{a^2}{\sigma ^2}\right)=1, \\
&c_0 \, _1F_1\left(-\frac{i \alpha }{2 \theta },\frac{1}{2},\frac{b^2}{\sigma ^2}\right)+c_1 b \, _1F_1\left(\frac{\theta -i \alpha }{2 \theta }, \frac{3}{2},\frac{b^2}{\sigma ^2}\right)=1.
\end{align}
\end{subequations}

Using these coefficients in \eqref{eq:gysol} and evaluating $g(x_0)$ results in the following for the characteristic function
\begin{subequations}\label{eq:psi}
\begin{equation}
\psi_{\tau}(\alpha)=\frac{N_{\psi}}{D_{\psi}},
\end{equation}
where
\begin{align}
N_{\psi}&=-x_0 \, _1F_1\left(-\frac{i \alpha }{2 \theta },\frac{1}{2},\frac{a^2}{\sigma ^2}\right) \, _1F_1\left(\frac{\theta -i \alpha }{2 \theta },\frac{3}{2},\frac{x_0^2}{\sigma ^2}\right) +a \, _1F_1\left(\frac{\theta -i \alpha }{2 \theta },\frac{3}{2},\frac{a^2}{\sigma ^2}\right) \, _1F_1\left(-\frac{i \alpha }{2 \theta },\frac{1}{2},\frac{x_0^2}{\sigma ^2}\right) \nonumber \\
& \qquad -b \, _1F_1\left(\frac{\theta -i \alpha }{2 \theta },\frac{3}{2},\frac{b^2}{\sigma ^2}\right) \, _1F_1\left(-\frac{i \alpha }{2 \theta },\frac{1}{2},\frac{x_0^2}{\sigma ^2}\right) +x_0 \, _1F_1\left(-\frac{i \alpha }{2 \theta },\frac{1}{2},\frac{b^2}{\sigma ^2}\right) \, _1F_1\left(\frac{\theta -i \alpha }{2 \theta },\frac{3}{2},\frac{x_0^2}{\sigma ^2}\right),\\
D_{\psi}&= a \, _1F_1\left(\frac{\theta -i \alpha }{2 \theta },\frac{3}{2},\frac{a^2}{\sigma ^2}\right) \, _1F_1\left(-\frac{i \alpha }{2 \theta },\frac{1}{2},\frac{b^2}{\sigma ^2}\right) -b \, _1F_1\left(-\frac{i \alpha }{2 \theta },\frac{1}{2},\frac{a^2}{\sigma ^2}\right) \, _1F_1\left(\frac{\theta -i \alpha }{2 \theta },\frac{3}{2},\frac{b^2}{\sigma ^2}\right).
\end{align}
\end{subequations}
The hypergeometric functions $_1F_1$ can be converted to other special functions, such as Hermite functions, parabolic cylinder functions, etc. \cite{abramowitz1964handbook}. Results on FPT of OU are presented in some of these forms  in standard texts \cite{BorodinSalminen12}.

%\subsection{Single threshold case}
\begin{remark}
If the FPT characteristic function is desired for a single threshold, it could be computed as a special case of the two threshold case analyzed here. There are two possibilities: either the initial condition is above the threshold, or below it. If the initial condition is above the threshold, then we may analyze this case as two thresholds case by letting the threshold at $b \to +\infty$, and considering $a$ as our threshold of interest. In the other case when the initial condition $x_0$ is below the threshold, then 
%we may analyze this as a special case of the two threshold setup. Here, 
we let $a \to -\infty$ and assume the threshold of interest at $b$.
\end{remark}

\begin{remark}
Recall our definition of the FPT for two threshold case given in \eqref{eq:fpt2thresholds}. The initial condition $x_0$ there is assumed to lie between the thresholds $a$ and $b$. If that were not the case, then the thresholds problem also becomes a single threshold problem. More specifically, if $x_0<a<b$, then the process will always reach $a$ before $b$. Therefore, the FPT is same as that for a single threshold at $a$. Analogously, if $x_0>b>a$, then the process will hit the threshold $b$ before the threshold $a$, and the FPT is same as that for reaching a single threshold at $b$.
\end{remark}

\subsection{When one of the threshold is reflecting}
Another possible situation of interest arises when one of the thresholds is reflecting. For example, we could assume that the threshold at $a$ is not absorbing and the process is reflected back as soon as it hits $a$. We are interested in computing the characteristic function of the first time at which the process reaches the threshold $b$. 

The computation follows the same principles as those for the two threshold case, and therefore reduces to solving the differential equation \eqref{eq:gODE} for $g(y)$. The general form of the solution in \eqref{eq:gysol} can be used in this case, with appropriate boundary conditions given by $g(b)=1$ and $g'(a)=0$ (see \cite{CoxMiller77,gardiner2004handbook} for more details). Note that if the absorbing threshold is at $a$ and $b$ is the reflecting threshold, then we will have the boundary conditions $g(a)=1$ and $g'(b)=0$. We do not analyze this case here.

Let us denote the FPT characteristic function $\psi_{\tau_r}$. Using the initial conditions to compute $c_0$ and $c_1$ in \eqref{eq:gysol} and then evaluating $g(x_0)$ results in the following for the characteristic function
\begin{subequations}\label{eq:psiRef}
\begin{equation}
\psi_{\tau r}(\alpha)=\frac{N_{\psi r}}{D_{\psi r}},
\end{equation}
where
\begin{align}
N_{\psi r}&= \frac{6 i \alpha}{\theta}  b x_0 \, _1F_1\left(1-\frac{i \alpha }{2 \theta },\frac{3}{2},\frac{a^2}{\sigma ^2}\right) \, _1F_1\left(\frac{\theta -i \alpha }{2 \theta },\frac{3}{2},\frac{x_0^2}{\sigma ^2}\right) \nonumber \\
&+2 a^2 \left(1-\frac{i\alpha}{\theta}\right) \, _1F_1\left(\frac{3}{2}-\frac{i \alpha }{2 \theta },\frac{5}{2},\frac{a^2}{\sigma ^2}\right) \, _1F_1\left(-\frac{i \alpha }{2 \theta },\frac{1}{2},\frac{x_0^2}{\sigma ^2}\right)+3 \sigma ^2 \, _1F_1\left(\frac{\theta -i \alpha }{2 \theta },\frac{3}{2},\frac{a^2}{\sigma ^2}\right) \, _1F_1\left(-\frac{i \alpha }{2 \theta },\frac{1}{2},\frac{x_0^2}{\sigma ^2}\right),
\end{align}
\begin{align}
D_{\psi r}&= 2 a^2 \left(1-\frac{i\alpha}{\theta}\right)\, _1F_1\left(\frac{3}{2}-\frac{i \alpha }{2 \theta },\frac{5}{2},\frac{a^2}{\sigma ^2}\right) \, _1F_1\left(-\frac{i \alpha }{2 \theta },\frac{1}{2},\frac{b^2}{\sigma ^2}\right) \nonumber \\
&+6 \frac{i\alpha}{\theta}  b^2 \, _1F_1\left(1-\frac{i \alpha }{2 \theta },\frac{3}{2},\frac{a^2}{\sigma ^2}\right) \, _1F_1\left(\frac{\theta -i \alpha }{2 \theta },\frac{3}{2},\frac{b^2}{\sigma ^2}\right)+3 \sigma ^2 \, _1F_1\left(\frac{\theta -i \alpha }{2 \theta },\frac{3}{2},\frac{a^2}{\sigma ^2}\right) \, _1F_1\left(-\frac{i \alpha }{2 \theta },\frac{1}{2},\frac{b^2}{\sigma ^2}\right).
\end{align}
\end{subequations}

So far we have computed the characteristic functions for FPT of OU process in various scenarios. The characteristic function can now be used to explore how various parameters affect the FPT statistics, and how they could be tuned to achieve desired FPT behavior. 

\begin{figure}[h]
\centering
\includegraphics[width=0.8\textwidth]{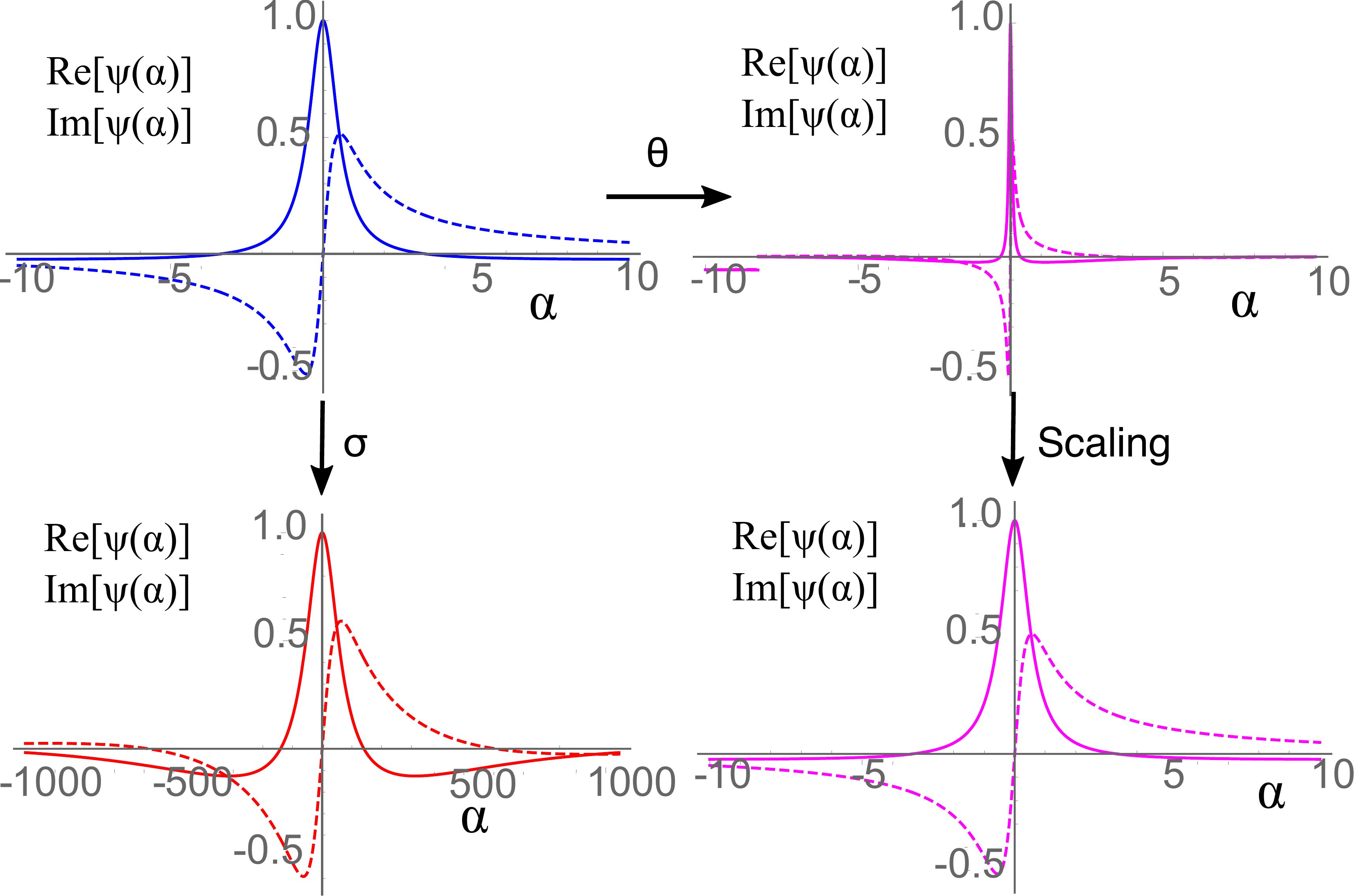}
\caption{Effect of the parameters $\theta$ and $\sigma$ on the shape of the characteristic function of the FPT for the case when both thresholds are absorbing. Both the real part (continuous line) and imaginary parts (dashed line) are shown. Upon changing the parameter $\theta$, the shape of the distribution does not change, and only its scale with respect to $\alpha$ changes. Upon scaling $\alpha$ by the same $\theta$, the previous shape can be reproduced. Upon scaling with the parameter $\sigma$, the shape of the distribution changes (note the lobe of the real part).}
\label{fig:scaleInv}
\end{figure}

\section{OPTIMAL PARAMETERS FOR DESIRED FPT MOMENTS}\label{sec:OptMom}
In this section, we investigate the effect of various parameters of the OU process on the FPT moments. Then, we examine how the parameters could be tuned so as to get a desired FPT moments.

\subsection{Scale invariance of the FPT distribution}

In the previous section, we derived characteristic functions of FPT distribution of the OU process under different scenarios (both thresholds absorbing or one of them reflecting). More generally, the characteristic function for other scenarios can also be derived from the generalized form in \eqref{eq:gysol}, with appropriate boundary conditions.  An important point to note is that in both \eqref{eq:psi} and \eqref{eq:psiRef}, the drift parameter $\theta$ always appears with $\alpha$ (as in $\alpha/\theta$). Therefore, if we consider the rescaled variable $\overline{\tau}=\tau \theta$, and find a general form similar to \eqref{eq:gysol}, it would be given by
\begin{subequations}
\begin{align}
\overline{g}(y)&=\E \left[ e^{i \alpha \overline{\tau} (y)} \right]=\E \left[ e^{i \alpha \theta \tau (y)} \right] \\
&=c_0 \, _1F_1\left(-\frac{i \alpha }{2},\frac{1}{2},\frac{y^2}{\sigma ^2}\right)+c_1 y \, _1F_1\left(\frac{1-i \alpha }{2}, \frac{3}{2},\frac{y^2}{\sigma ^2}\right).
\end{align}
\end{subequations}

Thus, the general solution $\overline{g}(y)$ for the rescaled variable $\overline{\tau}=\theta \tau$ would not depend on $\theta$. As the coefficients $c_0$ and $c_1$ above are obtained from boundary conditions, they would also be independent of $\theta$. 

An alternate way to infer this feature is to look at \eqref{eq:ou}. As $dw$ is of the order of $\sqrt{dt}$, we can rescale time by $\theta$ (as in $\overline{t}=\theta t$) and rewrite \eqref{eq:ou} as
\begin{equation}
dx=-x d\overline{t} + \sigma dw_{\overline{t}}.
\end{equation}
Because $\theta$ does not appear in the new time scale, the characteristic function of the FPT with this rescaling should be independent of $\theta$ as well.

To understand the implications of this property, consider the characteristic function of the rescaled variable $\overline{\tau}$
\begin{equation}
\psi_{\overline{\tau}}(\alpha)=\E \left[ e^ {i \alpha \overline{\tau}}\right]=1+i \alpha \E \left[\overline{\tau}\right] + \frac{i^2\alpha^2}{2!}\E \left[\overline{\tau^2}\right]+\ldots.
\end{equation}

Since $\theta$ does not appear in the above characteristic function, all moments of $\overline{\tau}$ are independent of $\theta$. Furthermore, because $\overline{\tau}=\theta \tau$, we have that
\begin{equation}
\E[\tau^m]=\theta^m \E[\overline{\tau}^m], \quad m \geq 1.
\end{equation}
Since $\E[\overline{\tau}]$ does not depend upon upon $\theta$, this implies that
\begin{equation}
\E[\tau^m] \propto \theta^m,
\end{equation}
and appropriately scaled moments of the FPT, $\E[\tau^m]/(\E[\tau])^m$, are independent of $\theta$. It follows that if we operate with normalized higher statistical moments such as the coefficient of variation (CV), skewness, kurtosis, etc, then changing the drift parameter $\theta$ only changes the mean FPT $\E[\tau]$.
%Because other moments scale accordingly with $\theta$, normalized statistical quantities such as the coefficient of variation (CV), skewness, kurtosis, etc. do not depend upon $\theta$. 
The scale invariance has been observed in distributions of other quantities \cite{giometto2013scaling}, and also of FPTs in other contexts \cite{redner2001guide,iyer2016first,ghusinga2016mechanistic}.

In terms of the characteristic function, we illustrate the scale invariance property in Fig.~\ref{fig:scaleInv} for the case when both thresholds are absorbing. The real and imaginary parts of the characteristic function are plotted for the FPT. By varying the drift parameter $\theta$, the characteristic function $\psi_{\tau}(\alpha)$ does not change in shape and just scales with respect to the $\alpha$ axis. However, changing the relative noise strength $\sigma$ affects the shape of the characteristic function. Similar behavior is also seen in the case when one of the thresholds is reflecting, though the results are not shown in order to avoid repetition.

\subsection{Tuning FPT moments}
Recall the form of \eqref{eq:ou}. Suppose that we are interested in tuning the two parameters ($\theta$ and $\sigma$) of the process so as to get desired moments of the FPT. Since $\theta$ only changes the mean, and the other quantities of interest (such as coefficient of variation (CV), skewness etc.) are independent of it, one could independently tune the mean FPT and one other quantity. Typically, the CV is the other quantity of interest because it represents the noise in the FPT.

What remains to be seen is how the relative noise strength $\sigma$ affects the mean and CV of the FPT. One could then choose appropriate $\sigma$ such that the CV is at a desired level, and then tune $\theta$ to get the desired mean. It turns out that both the mean and CV are decreasing functions of $\sigma$ for the two absorbing thresholds case. The CV eventually approaches to a limiting value
\begin{equation}
\sqrt{\frac{(x_0-a)^2+(b-x_0)^2}{3(x_0-a)(b-x_0)}},
\end{equation}
 which corresponds to the CV of the FPT for a diffusion with zero drift. In case when the threshold at $a$ is reflecting, the mean still decreases with increase in $\sigma$. The CV, on the other hand, shows a slight dip before increasing to a limiting value
\begin{equation}
\sqrt{\frac{2 \left( (x_0-a)^2+(b-a)^2\right)}{3 (b-x_0) \left((x_0-a)+(b-a)\right)}}
\end{equation}
that corresponds to the CV of the FPT for a diffusion with zero drift.

\begin{figure}[h]
\centering
\includegraphics[width=0.5\linewidth]{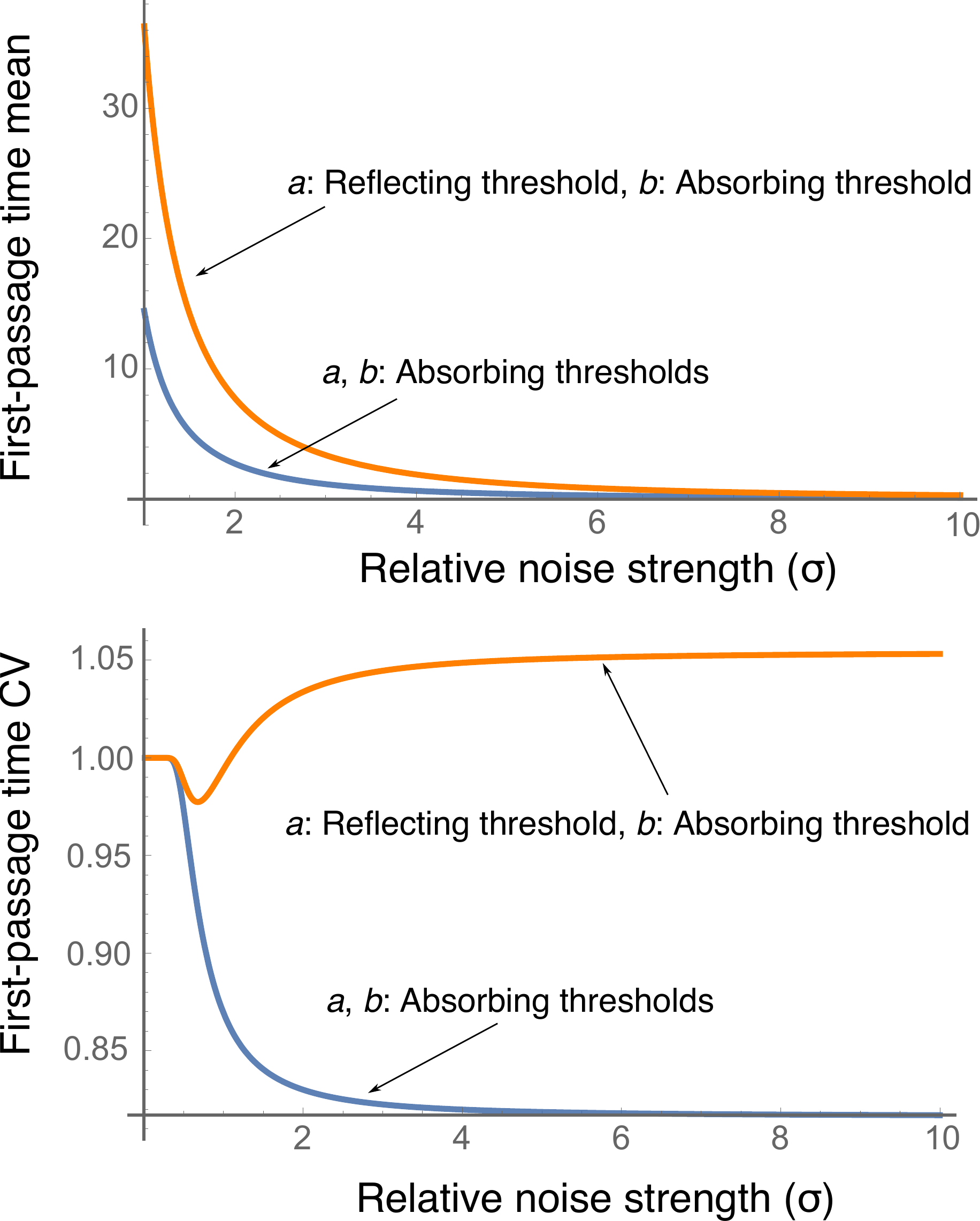}
\caption{Effect of the relative noise strength ($\sigma$) on the mean and coefficient of variation of FPT.}
\label{fig:meanCVwithSigma}
\end{figure}

Collectively, these results show that if one were to tune $\theta$, and $\sigma$ then any desired mean FPT could be achieved, but there is a limit to the achievable CV. Achieving a low CV in the double thresholds case requires a high value of $\sigma$ whereas for the case when one the barriers is reflecting, there is an optimal $\sigma$ that minimizes the CV.

\subsection{Effect of thresholds and initial condition}
\begin{figure}[h]
\centering
\includegraphics[width=0.5\linewidth]{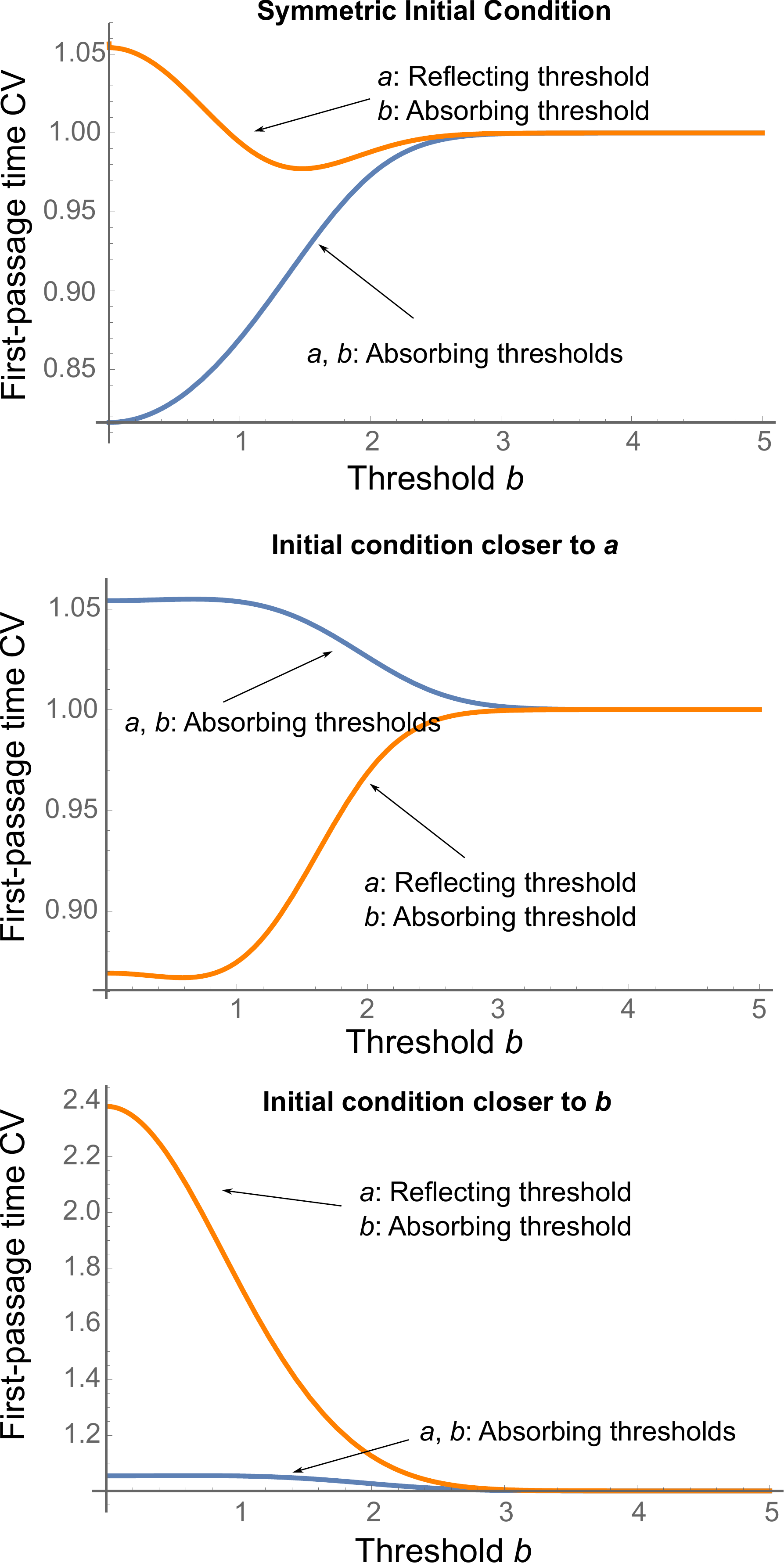}
\caption{Effect of varying threshold on coefficient of variation for symmetric and asymmetric initial conditions.}
\label{fig:CVwithz}
\end{figure}

Our analysis thus far has assumed fixed thresholds and a given initial condition. In Fig.~\ref{fig:CVwithz}, we examine how the results change when one of these parameters are changed. As a first case, we consider a symmetric thresholds, i.e., $a=-b$, and the initial condition to be at $x_0=0$. In this case, increasing $b$ leads to increase in CV of FPT if both thresholds are absorbing. In contrast, if $a$ is considered to be reflecting, then there is an optimal threshold $b$ at which the CV hits a minimum.
%before increasing back again. 
Increasing the threshold beyond a certain point does not affect the CV anymore. This corresponds to the situation when the absorbing threshold(s) is far from the initial condition and crossing it is dominated purely by noise (see Fig.~\ref{fig:CVwithz}, top).

Next, we consider the case when the initial condition $x_0$ is not symmetric. Assuming the threshold to be $a=-b$, we take two cases: $x_0=-\frac{b}{2}$ and $x_0=\frac{b}{2}$. When both $a$ and $b$ are absorbing, increasing the threshold decreases CV of FPT and the CV seems to approach the limit of symmetric initial condition (Fig.~\ref{fig:CVwithz}, middle and bottom). However, when the threshold $a$ is taken as reflecting, then the CV properties change depending upon $x_0$. More specifically, when $x_0$ is near the reflecting threshold, then increasing the threshold increases CV. In contrast, when $x_0$ is near the absorbing threshold, then increasing the threshold leads to reduction in CV of FPT.

To sum up, the FPT distribution for an OU process is scale invariant with respect to the drift parameter, and thereby allows independent tuning of the mean FPT and another statistical quantity that consists of appropriately scaled moments of the FPT. If one is interested in obtaining a FPT distribution that matches more than two statistical quantities of interest, it is not possible. A question of interest at this point is how close can the FPT distribution get to a given distribution?

\section{OPTIMAL PARAMETERS FOR DESIRED FPT DISTRIBUTION}\label{sec:OptDist}
Suppose that instead of tuning the moments, we are interested in tuning the distribution of the FPT itself. More specifically, we are interested in choosing the parameters such that the FPT distribution is as close to a desired distribution as possible. In this section, we discuss the tuning of OU process to achieve such behavior. 

To this end, we consider the relation between probability density function and the characteristic function stated in \eqref{eq:distmetric}. Although the desired distribution could be specified as any distribution of interest, we consider the Gamma distribution here. The rationale behind this is that the Gamma distribution is
that it is the distribution of a summation of exponential random variables and in a limiting case, it can even represent a degenerate (deterministic) distribution.

\begin{figure}[h]
\centering
\includegraphics[width=0.5\linewidth]{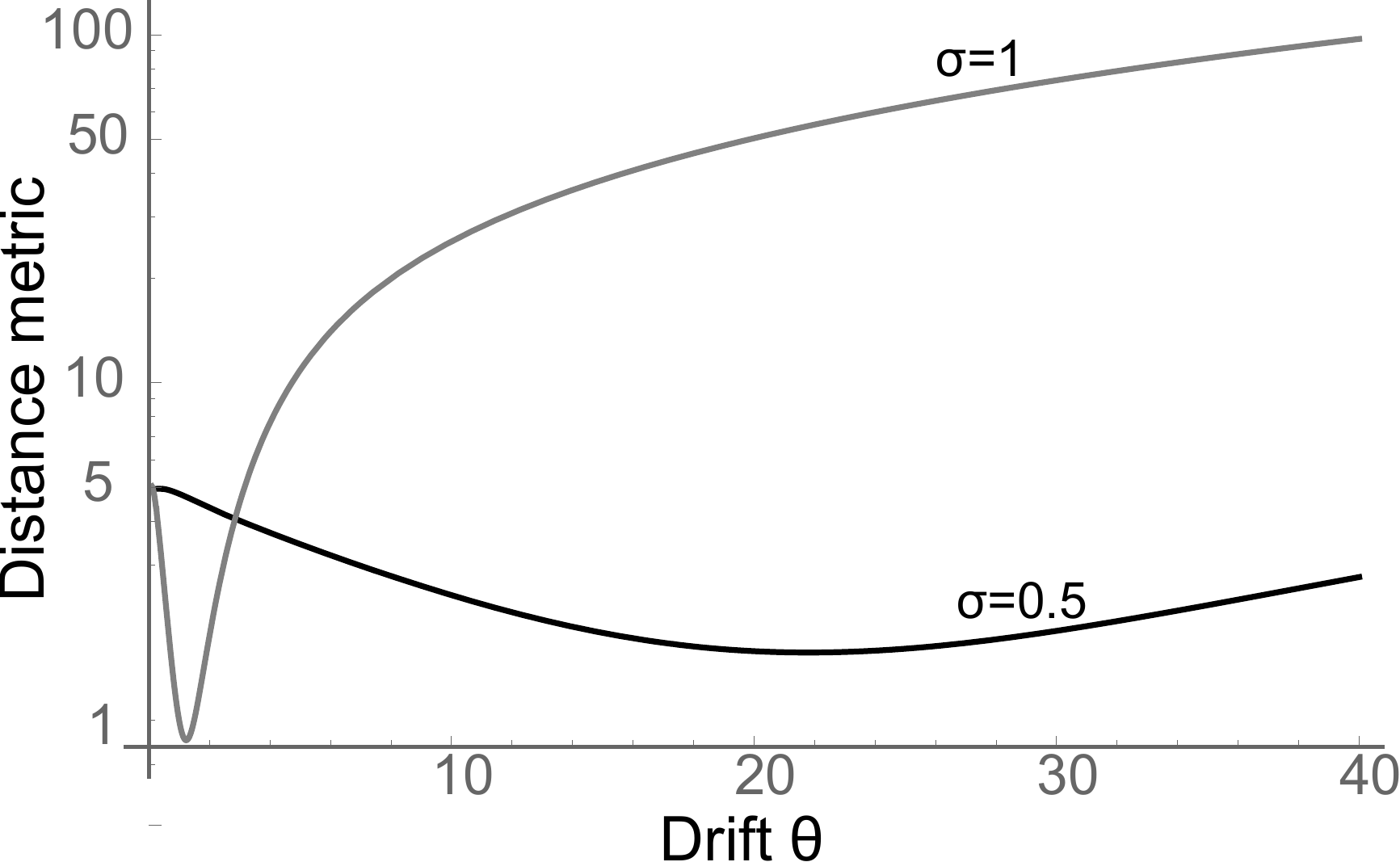}
\caption{The distance between the FPT distribution and the desired distribution is minimized at an optimal drift parameter $\theta$ for a given relative noise strength $\sigma$. Increasing $\sigma$ results in lower value of $\theta$ at which the distance is minimized, and the distance metric itself also reduces in this case.}
\label{fig:distanceOpt}
\end{figure}

In Fig.~\ref{fig:distanceOpt}, we assume both thresholds to be absorbing and plot the distance metric between the desired distribution and the FPT distribution as a function of $\theta$. The value of $\sigma$ is taken to be fixed. It can be seen that there is an optimal value of $\theta$ that minimizes the distance metric. Furthermore, if $\sigma$ is increased and this process is iterated, we see that the optimal value of $\theta$ decreases, and the minimum of the distance metric decreases as well. Referring to Fig.~\ref{fig:scaleInv}, it can be seen that fixing the value of $\sigma$ basically fixes the shape of the characteristic function, and then by changing $\theta$, an appropriate scale is chosen so that the characteristic function matches with the desired one. By iterating over $\sigma$, we change the shape of the characteristic function and find the optimal shape that matches the desired characteristic function with appropriate scaling $\theta$. Even in the case when one of the boundary is reflecting, we get optimal parameters that minimize the distance metric. The results are not presented here to avoid redundancy.

% A plot of the characteristic functions for the desired distribution and the first-passage time distribution is shown in Fig.~\ref{}. Clearly, they do not match perfectly, as expected because only the parameter $\sigma$ changes the shape of the distribution, and $\theta$ just shifts the mean to the desired mean level. Although the means do not match exactly as there is a trade off in terms of distance from the other moments.

% \com{need to put a plot showing minimum of the distance metric is achieved at some value of $\sigma$ and $\theta$. should consider both types of thresholds, and symmetric and asymmetric initial conditions.}
%\subsection{Single threshold}
%\subsection{Absorbing and reflecting barriers}

\section{CONCLUSION}\label{sec:conclusion}
The OU process is used to model stochastic phenomena in a variety of contexts. In particular, the FPT of OU process has been used to study decision-making~\cite{bogacz2006physics}, animal foraging~\cite{OB-MC-MM-PHS-RV:05}, financial markets~\cite{Lin1998}, and synchronization of clocks~\cite{WS-MZW-LAS:12}. In this paper, we consider OU process with two thresholds. Both of the thresholds could be absorbing, or one of them could be reflecting while the other one is absorbing. We analyze the effect of OU process parameters on its FPT statistics, and how the parameters could be chosen to obtain a desired FPT statistics. 

Future work will focus on analyzing the optimal parameters by not restricting $\theta$ to be positive, and thereby allowing the OU process to not just be mean reverting. Given the numerical tractability of characteristic functions, it would also be interesting to explore other stochastic processes the approach presented here and explore optimal/sub-optimal control strategies that result in desired FPT statistics.

%\section*{Appendix}

\begin{scratch}
%\section{Notes}

\end{scratch}

% The mean first-passage time can be computed as the solution to the following differential equation \cite{GihmanSkorohodSDE72}
% \begin{equation}
% \frac{1}{2}\sigma^2 \ddot{v}(y)+\theta(\mu-y)\dot{v}(y)=-1,
% \end{equation}
% with $v(a)=v(b)=0$. Furthermore, the second order moment can be computed as the solution to the differential equation \cite{GihmanSkorohodSDE72}
% \begin{equation}
% \frac{1}{2}\sigma^2 \ddot{\eta}(y)+\theta(\mu-y)\dot{\eta}(y)=-2v(y),
% \end{equation}
% with $\eta(a)=\eta(b)=0$.
%\section*{Acknowledgment}
%The authors thank...

%\section*{References}
\bibliographystyle{IEEEtran}
\bibliography{bibLoc}

\begin{thebibliography}{10}
\providecommand{\url}[1]{#1}
\csname url@samestyle\endcsname
\providecommand{\newblock}{\relax}
\providecommand{\bibinfo}[2]{#2}
\providecommand{\BIBentrySTDinterwordspacing}{\spaceskip=0pt\relax}
\providecommand{\BIBentryALTinterwordstretchfactor}{4}
\providecommand{\BIBentryALTinterwordspacing}{\spaceskip=\fontdimen2\font plus
\BIBentryALTinterwordstretchfactor\fontdimen3\font minus
  \fontdimen4\font\relax}
\providecommand{\BIBforeignlanguage}[2]{{%
\expandafter\ifx\csname l@#1\endcsname\relax
\typeout{** WARNING: IEEEtran.bst: No hyphenation pattern has been}%
\typeout{** loaded for the language `#1'. Using the pattern for}%
\typeout{** the default language instead.}%
\else
\language=\csname l@#1\endcsname
\fi
#2}}
\providecommand{\BIBdecl}{\relax}
\BIBdecl

\bibitem{bogacz2006physics}
R.~Bogacz, E.~Brown, J.~Moehlis, P.~Holmes, and J.~D. Cohen, ``The physics of
  optimal decision making: a formal analysis of models of performance in
  two-alternative forced-choice tasks.'' \emph{Psychological review}, vol. 113,
  p. 700, 2006.

\bibitem{HCT-FYMW:00}
H.~C. Tuckwell and F.~Y.~M. Wan, ``First passage time to detection in
  stochastic population dynamical models for {HIV-1},'' \emph{Applied
  Mathematics Letters}, vol.~13, no.~5, pp. 79--83, 2000.

\bibitem{GV-BF-JC-AD-ERD:08}
G.~Vahedi, B.~Faryabi, J.~Chamberland, A.~Datta, and E.~R. Dougherty,
  ``Intervention in gene regulatory networks via a stationary
  mean-first-passage-time control policy,'' \emph{IEEE Transactions on
  Biomedical Engineering}, vol.~55, no.~10, pp. 2319--2331, 2008.

\bibitem{SS-DJS-CRM:14}
S.~Singh, D.~J. Schneider, and C.~R. Myers, ``Using multitype branching
  processes to quantify statistics of disease outbreaks in zoonotic
  epidemics,'' \emph{Physical Review E}, vol.~89, no.~3, p. 032702, 2014.

\bibitem{redner2001guide}
S.~Redner, \emph{A guide to first-passage processes}.\hskip 1em plus 0.5em
  minus 0.4em\relax Cambridge University Press, 2001.

\bibitem{OB-MC-MM-PHS-RV:05}
O.~B{\'e}nichou, M.~Coppey, M.~Moreau, P.~Suet, and R.~Voituriez, ``Optimal
  search strategies for hidden targets,'' \emph{Physical Review Letters},
  vol.~94, no.~19, p. 198101, 2005.

\bibitem{WS-MZW-LAS:12}
W.~Suwansantisuk, M.~Z. Win, and L.~A. Shepp, ``First passage time problems
  with applications to synchronization,'' in \emph{IEEE International
  Conference on Communications}, 2012, pp. 2580--2584.

\bibitem{Farkas2001}
Z.~Farkas and T.~Fulop, ``One-dimensional drift-diffusion between two absorbing
  boundaries: application to granular segregation,'' \emph{Journal of Physics
  A: Mathematical and General}, vol.~34, no.~15, pp. 3191--3198, 2001.

\bibitem{Lin1998}
X.~S. Lin, ``Double barrier hitting time distributions with applications to
  exotic options,'' \emph{Insurance: Mathematics and Economics}, vol.~23,
  no.~1, pp. 45--58, 1998.

\bibitem{Horrocks2004}
J.~Horrocks and M.~E. Thompson, ``Modeling event times with multiple outcomes
  using the {Wiener} process with drift,'' \emph{Lifetime Data Analysis},
  vol.~10, no.~1, pp. 29--49, 2004.

\bibitem{chen2016optimal}
Y.~Chen, T.~T. Georgiou, and M.~Pavon, ``Optimal steering of a linear
  stochastic system to a final probability distribution, part i,'' \emph{IEEE
  Transactions on Automatic Control}, vol.~61, no.~5, pp. 1158--1169, 2016.

\bibitem{gardiner2004handbook}
C.~Gardiner, ``Handbook of stochastic methods for physics, chemistry, and the
  natural sciences,'' \emph{Springer series in synergetics}, 2004.

\bibitem{sacerdote2003threshold}
L.~Sacerdote and C.~Zucca, ``Threshold shape corresponding to a {G}amma firing
  distribution in an ornstein-uhlenbeck neuronal model,'' \emph{Scientiae
  Mathematicae Japonicae}, vol.~58, no.~2, pp. 295--306, 2003.

\bibitem{krstic2008boundary}
M.~Krstic and A.~Smyshlyaev, \emph{Boundary control of PDEs: A course on
  backstepping designs}.\hskip 1em plus 0.5em minus 0.4em\relax SIAM, 2008,
  vol.~16.

\bibitem{Ghusinga17}
K.~R. Ghusinga, J.~J. Dennehy, and A.~Singh, ``First-passage time approach to
  controlling noise in the timing of intracellular events,'' \emph{Proceedings
  of the National Academy of Sciences}, vol. 114, pp. 693--698, 2017.

\bibitem{BorodinSalminen12}
A.~N. Borodin and P.~Salminen, \emph{Handbook of Brownian motion-facts and
  formulae}.\hskip 1em plus 0.5em minus 0.4em\relax Birkh{\"a}user, 2012.

\bibitem{UshakovCF99}
N.~G. Ushakov, \emph{Selected topics in characteristic functions}.\hskip 1em
  plus 0.5em minus 0.4em\relax Walter de Gruyter, 1999.

\bibitem{GihmanSkorohodSDE72}
I.~I. Gihman and A.~V. Skorohod, \emph{Stochastic Differential
  Equations}.\hskip 1em plus 0.5em minus 0.4em\relax Springer-Verlag Berlin
  Heidelberg, 1972.

\bibitem{CoxMiller77}
D.~R. Cox and H.~D. Miller, \emph{The theory of stochastic processes}.\hskip
  1em plus 0.5em minus 0.4em\relax CRC Press, 1977.

\bibitem{abramowitz1964handbook}
M.~Abramowitz and I.~A. Stegun, \emph{Handbook of mathematical functions: with
  formulas, graphs, and mathematical tables}.\hskip 1em plus 0.5em minus
  0.4em\relax Courier Corporation, 1964, vol.~55.

\bibitem{giometto2013scaling}
A.~Giometto, F.~Altermatt, F.~Carrara, A.~Maritan, and A.~Rinaldo, ``Scaling
  body size fluctuations,'' \emph{Proceedings of the National Academy of
  Sciences}, vol. 110, no.~12, pp. 4646--4650, 2013.

\bibitem{iyer2016first}
S.~Iyer-Biswas and A.~Zilman, ``First-passage processes in cellular biology,''
  \emph{Advances in Chemical Physics, Volume 160}, pp. 261--306, 2016.

\bibitem{ghusinga2016mechanistic}
K.~R. Ghusinga, C.~A. Vargas-Garcia, and A.~Singh, ``A mechanistic stochastic
  framework for regulating bacterial cell division,'' \emph{Scientific
  Reports}, vol.~6, 2016.

\end{thebibliography}

\end{document}